\documentclass[11pt]{article}
\usepackage{latexsym,amsthm,verbatim,ifthen,graphicx,color,mathrsfs}
\usepackage[all]{xy}
\usepackage{color}
\usepackage[symbol]{footmisc}

\usepackage{color}

\usepackage{lscape}
 \usepackage[ansinew]{inputenc}
 \usepackage{multicol}
 \usepackage[all]{xy}\xyoption{poly}
 \usepackage{mathrsfs}
\usepackage[psamsfonts]{eucal}
\usepackage{amssymb, amsmath}
\usepackage{geometry}
\usepackage{enumerate}
\usepackage{float}

\newtheorem{theorem}{Theorem}[section]

\newtheorem{corollary}[theorem]{Corollary}
 \newtheorem{lemma}[theorem]{Lemma}
 \newtheorem{proposition}[theorem]{Proposition}
 \theoremstyle{definition}
 \newtheorem{definition}[theorem]{Definition}

 \newtheorem{rem}[theorem]{Remark}

\def\quil{{\mathscr L}}

\def\quic{{\mathscr C}}

\def\esp#1{{\quad\text{#1}\quad}}

\newcommand{\bz}{\mathbb Z}

\def\bq{\mathbb{Q}}

\def\des{{s^{-1}}}

\newcommand{\coker}{{\rm{coker}\,}}
    \newcommand{\lasu}{{\mathfrak{L}}}

 \newcommand{\lib }{\mathbb{L}}

  \newcommand{\coprodc}{\,\widehat{\amalg }\,\,}

  \newcommand{\Hom}{\operatorname{\text{\rm Hom}}}

\newcommand{\catcdgl}{\operatorname{{\bf cdgl}}}

\newcommand{\catss}{\operatorname{{\bf sset}}}

\newcommand{\Ho}{\operatorname{{\rm Ho}}}

 \newcommand{\map}{\operatorname{{\rm map}}}

    \newcommand{\id}{\operatorname{{\rm id}}}
        \newcommand{\ad}{\operatorname{{\rm ad}}}

 \newcommand{\MC}{\operatorname{{\rm MC}}}

\newcommand{\mc}{{\MC}}

  \newcommand{\otimesc}{\widehat{\otimes}}

   \newcommand{\libc}{{\widehat\lib}}

\usepackage{hyperref}
\usepackage{physics}
\hypersetup{
    colorlinks=true,
    linkcolor=blue,
    filecolor=magenta,
    urlcolor=cyan,
    citecolor=blue,
}

\newcommand{\acento}{{\scriptscriptstyle \wedge}}

\theoremstyle{remark}

\usepackage{stmaryrd}

\usepackage{ textcomp }
\newcommand{\cdgl}{\operatorname{\bf cdgl}}

\begin{document}

\title{A Lie characterization of the Bousfield-Kan $\bq$-completion and $\bq$-good spaces}

\author{Yves F\'elix, Mario Fuentes and Aniceto Murillo\footnote{The authors were partially supported by the Spanish Government grant PID2023-149804NB-I00. The second author was also supported by the Spanish Government grant JDC2024-053124-I.}}

\maketitle

\begin{abstract} The model and realization functors of Buijs-F\'elix-Murillo-Tanr\'e provide a Quillen adjunction $
 \xymatrix{\catss& \catcdgl \ar@<0.5ex>[l]^(.47){\langle\,\cdot\,\rangle}
\ar@<0.5ex>[l];[]^(.50){\lasu}\\}
$ between simplicial sets and complete differential graded Lie algebras. We prove that the unit of this adjunction is, up to homotopy, the Bousfield-Kan $\bq$-completion. As consequences, we obtain Lie-theoretic characterizations of $\bq$-good spaces and equivalences between suitable homotopy categories of rational spaces and complete differential graded Lie algebras.
\end{abstract}

\section*{Introduction}

The Bousfield-Kan $\bq$-completion \cite{bouskan} is a coaugmented functor
$$
(\cdot)^\acento_\bq\colon \catss\longrightarrow\catss
$$
which extends the classical rationalization for nilpotent spaces and it satisfies the following property: a map $f\colon X\to Y$ induces an isomorphism on rational homology if and only if $f^\acento_\bq\colon X^\acento_\bq\to Y^\acento_\bq$ is a weak homotopy equivalence.

On the other hand, in  \cite{bufemutan0}, the authors introduce a Quillen pair, model and realization,
 $$
 \xymatrix{\catss& \catcdgl \ar@<0.75ex>[l]^(.47){\langle\,\cdot\,\rangle}
\ar@<0.75ex>[l];[]^(.50){\lasu}\\}
$$
between the categories of simplicial sets and complete differential graded Lie algebras (cdgl's henceforth) whose restriction
$$
\xymatrix{ \catss_0& \catcdgl_0 \ar@<0.75ex>[l]^(.50){\langle\,\cdot\,\rangle}
\ar@<0.75ex>[l];[]^(.50){\lasu^a}\\}
$$
to reduced simplicial sets and connected cdgl's
extends, up to homotopy, the classical model and realization functors of Quillen \cite{qui}. Here, $\lasu^a$ is the component of $\lasu$ at the $0$-simplex $a$ (see \S1).

 We prove:

\begin{theorem}\label{main} For any reduced simplicial set $X$, the unit of the adjunction
$$
X\longrightarrow \langle\lasu^a_X\rangle
$$
is weakly equivalent to the Bousfield-Kan $\bq$-completion of $X$
$$
X\longrightarrow X^\acento_\bq.
$$
\end{theorem}
In general, for any  simplicial set $X=\amalg_{i\in I} X_i$ not necessarily connected,  $\langle \lasu_X\rangle$ is weakly homotopy equivalent to  $\amalg_{i\in I}\langle \lasu_{X_i}^{a_i}\rangle\coprod\{*\}$, with $a_i$ any vertex of $X_i$ for all $i\in I$. Thus, the previous result provides:

\begin{corollary}
For any simplicial set $X$,
$
\langle \lasu_X\rangle\simeq \amalg_{X_i\in\pi_0(X)}\, {(X_i)}^\acento_\bq\amalg\{*\}$.
\hfill$\square$
\end{corollary}

It should be noted that, in the case where $X$ is a reduced simplicial set of finite type, the identification in Theorem \ref{main} of its $\bq$-completion with the realization of its Lie model already appears in Theorem 11.14 of \cite{bufemutan0}. Indeed, this result follows from the identification of $\langle \lasu_X^a\rangle$ with the homotopy type of the geometric realization of the Sullivan cdga model of $X$, together with the key result, Theorem 12.2 of \cite{bousgu}, where Bousfield and Gugenheim prove that this realization is homotopy equivalent to $X^\acento_\bq$. As it is inherent in Sullivan's approach, finite type assumptions are essential.

 Theorem \ref{main} extends this result to any arbitrary reduced simplicial set $X$, providing a self-contained, Lie-theoretic description of the Bousfield-Kan $\bq$-tower of $X$ (see \cite[Chapter III, \S6]{bouskan}) whose inverse limit yields $X^\acento_\bq$. In our approach this tower is expressed explicitly  in terms of the lower central series quotients of the Lie model of $X$.

\medskip

On the other hand, recall that the $\bq$-completion is not an idempotent functor. In fact $X^\acento_\bq\simeq (X^\acento_\bq)^\acento_\bq$ if and only if the natural map $X\to X^\acento_\bq$ induces a rational homology isomorphism.  These are the so called {\em $\bq$-good} spaces. For instance, and quite remarkably, the wedge of two circles is a $\bq$-bad space \cite{IM}. It turns out that $X$ is $\bq$-good if and only if $X^\acento_\bq$ is $\bq$-good, and this happens if  and only if $X\to X^\acento_\bq$ is the Bousfield $\bq$-homology localization \cite{bous}. Hence, Theorem \ref{main} has the following direct consequence:
\begin{corollary}
(i) A  reduced simplicial set $X$ is $\bq$-good if and only if the unit $
X\to\langle\lasu_X^a\rangle
$ induces an isomorphism in rational homology.

(ii)
The model functor induces an injective-on-objects functor
 $$
 \xymatrix{ \Ho\catss^{\scriptscriptstyle\wedge}_{0,\bq} \ar[r]^(.50){\lasu^a}& \Ho\catcdgl_0
\\}
$$
between the homotopy categories of reduced,  $\bq$-good, $\bq$-local simplicial sets and connected cdgl's. \hfill$\square$
\end{corollary}
 Also, Theorem \ref{main} together with  \cite[Theorem~6.6]{ber}, see also \cite[Theorem~12.18]{bufemutan0}, immediately imply the following, under no finiteness restrictions:

 \begin{corollary} For any cdga model $A$ and any Lie model $L$ of the reduced simplicial sets $X$ and $Y$ respectively, the realization of the cdgl $A\otimesc L$ is weakly equivalent to the mapping space $\map(X,Y^\acento_\bq)$.\hfill$\square$
\end{corollary}

Finally, upon restricting to nilpotent spaces, we obtain the following particular consequence of the preceding results:

\begin{theorem}\label{main2}
The model and realization functors induce equivalences
 $$
 \xymatrix{ \Ho\catss^{\mathcal N}_{0,\bq}& \Ho\catcdgl^{\mathcal N}_0 \ar@<0.75ex>[l]^(.47){\langle\,\cdot\,\rangle}
\ar@<0.75ex>[l];[]^(.50){\lasu^a}\\}
$$
between the homotopy categories of reduced nilpotent rational simplicial sets  and homologically nilpotent connected cdgl's.
\end{theorem}

Similarly, Theorem \ref{main2} also has several precursors in the finite type setting. Indeed, in \cite{nei} (see also Chapters II and III of \cite{tan} and \cite[\S10.1]{bufemutan0}), Neisendorfer established that the homotopy category of connected differential graded Lie algebras with nilpotent homology of finite type is equivalent to the rational homotopy category of nilpotent spaces with finite type rational homology. In these works, the finiteness assumptions arise again from the need  to relate differential graded Lie algebra models with Sullivan's commutative differential graded algebra approach to rational homotopy theory, which breaks down if finite type restrictions are not imposed. In contrast, Theorem \ref{main2} is obtained here as a direct consequence of the Lie-theoretic self-contained characterization of the Bousfield-Kan $\bq$-completion
 provided by Theorem \ref{main}.

\smallskip

As for the organization of the paper, Section 1 contains a brief summary of the relevant results from the homotopy theory of cdgl's that will be used throughout the paper. Section 2 is almost entirely devoted to the proof of Theorem \ref{main}.  Given a reduced simplicial set $X$, we first carry out an exhaustive study of the sequence given by the quotients of   $\lasu_X^a$ by the terms of its lower central series, together with the description of suitable models for these quotients. We then show that the realization of this sequence, which is precisely $\langle \lasu_X^a\rangle$, provides a $\bq$-tower of $X$ in the sense of \cite[Chapter III, \S6]{bouskan} and thus, its inverse limit yields $X^\acento_\bq$.

\section{Preliminaries}\label{prelimi}

All that follows can be found in the detailed reference \cite{bufemutan0} or the original sources \cite{bufemutan1,bufemutan2}.

The ground field of any considered vector space will always be $\bq$. A {\em filtration} of a differential graded Lie algebra (dgl for short)  $L$, or $(L,d)$ if we want to specify the differential,  is a decreasing sequence of differential Lie ideals,
$$L=F^1\supset\dots \supset F^n\supset F^{n+1}\supset\cdots$$ such that $[F^p,F^q]\subset F^{p+q}$ for $p,q\geq 1$. The lower central series,
$$
L^1\supset\dots\supset L^n\supset L^{n+1}\supset\cdots
$$
with $L^{1}= L$ and $L^{n}= [L, L^{n-1}]$ for $n>1$, is a fundamental example.
A {\em complete differential graded Lie algebra}, cdgl henceforth,  is a dgl $L$ equipped with a  filtration $\{F^n\}_{n\ge 1}$ for which the  natural map
$$
L\stackrel{\cong}{\longrightarrow}\varprojlim_n L/F^n
$$
is a dgl isomorphism. We denote by $\catcdgl$ the  category of cdgl's and morphisms which are dgl maps preserving the corresponding filtrations.

Given a  dgl  $L$  filtered by $\{F^n\}_{n\ge 1}$, its {\em completion} is the dgl
$
\widehat L=\varprojlim_nL/F^{n}
$
which is always complete with respect to the filtration
$
\widehat{F}^n=\ker ( \widehat L \to L/F^n)
$
since $ \widehat L/\widehat{F}^n=L/F^n$. In particular, denote by
$$
\libc(V)=\varprojlim_n\lib(V)/\lib(V)^n$$
 the completion of the free Lie algebra $\lib(V)$ generated by the graded vector space $V$.

\medskip

\noindent{\em Caution:} In the following it is important not to confuse the filtration degree with the length of a word. If we denote by $\libc^q(V)$ the vector spaces generated by words of bracket length exactly $q$, then
 $$
 \widehat{\mathbb L}^{\geq n} (V)=  \prod_{q\geq n} \mathbb L^q(V),\esp{for} n\ge 1.
 $$

A {\em Maurer-Cartan} or $\mc$ element of a given dgl $L$ is an element $a\in L_{-1}$ for which $da=-\frac{1}{2}[a,a]$.

There is a pair of adjoint functors, {\em (global) model} and {\em realization},
\begin{equation}\label{pair}
\xymatrix{ \catss& \catcdgl \ar@<0.75ex>[l]^(.50){\langle\,\cdot\,\rangle}
\ar@<0.75ex>[l];[]^(.50){\lasu}\\},
\end{equation}
which restricts to adjoint functors
\begin{equation}\label{pair2}
\xymatrix{ \catss_0& \catcdgl_0 \ar@<0.75ex>[l]^(.50){\langle\,\cdot\,\rangle}
\ar@<0.75ex>[l];[]^(.50){\lasu^a}\\}
\end{equation}
between the categories of reduced simplicial sets and connected cdgl's. Here, for any reduced simplicial set $X$,  $\lasu^a_X=\lasu_X/(a)$ where $(a)$ denotes the cdgl ideal generated by the only $0$-simplex  $a$ of $X$ which is naturally identified to a Maurer-Cartan element of $\lasu_X$.

These functors are
 based on the cosimplicial cdgl
$\lasu_\bullet=\{\lasu_n\}_{n\ge 0}$ which we briefly recall.
 For each $n\ge 0$, see \cite[Chapter 6]{bufemutan0},
$$
\lasu_n=\bigl(\libc(s^{-1}\Delta^{n}),d)
$$
in which
$s^{-1}\Delta^{n}$  denotes the desuspension $s^{-1}N_*(\underline\Delta^{n})$ of the non-degenerate simplicial chains  on $\underline\Delta^n$. Hence, generators of  degree $p-1$ of $s^{-1}\Delta^n$, with $p\ge 0$, are written as $a_{i_0\dots i_p}$ with $0\le i_0<\dots<i_p\le n$. The cofaces and codegeneracies in $\lasu_\bullet$ are induced by those on the cosimplicial chain complex $\des N_*\underline\Delta^\bullet$, and  the differential $d$ on each $\lasu_n$ is the only one (up to cdgl isomorphism) satisfying:
\begin{itemize}
 \item[(1)] For each $i=0,\dots,n$, the generators of $s^{-1}\Delta^{n}$, corresponding to vertices, are MC elements.

     \item[(2)] The linear part of $d$ is induced by the boundary operator of $s^{-1}\Delta^{n}$.
     \item[(3)] The cofaces and codegeneracies are cdgl morphisms.
     \end{itemize}

The realization of a given cdgl $L$ is the simplicial set
     $$
     \langle L\rangle=\Hom_{\catcdgl}(\lasu_\bullet,L).
     $$
 In particular, if $L$ is connected, $\langle L\rangle$ is  reduced and there are isomorphisms
 \begin{equation}\label{homoto}
 \pi_n\langle L\rangle\cong  H_{n-1}(L),\quad\text{for any $n\ge 1$.}
  \end{equation}
  Here, the group structure in $H_0(L)$ is considered with the Baker-Campbell-Hausdorff (BCH) product. Moreover, the action of $\pi_1\langle L\rangle$ on $\pi_*\langle L\rangle$ is identified to the {\em exponential action}
 \begin{equation}\label{action}
 H_0(L)\times H(L)\longrightarrow H(L),\quad (\alpha,\beta)\mapsto e^{\ad_\alpha}(\beta).
 \end{equation}

     On the other hand, the global model of a simplicial set $X$ is the cdgl
     $$
     \lasu_X=\varinjlim_{\sigma\in X}\lasu_{|\sigma|}.
     $$
As a complete Lie algebra,   $
\lasu_X=\libc(s^{-1}X)
$ where  $s^{-1}X$ denotes the desuspension $s^{-1}N_*(X)$ of the chain complex of non-degenerate simplicial chains on $X$. Moreover, the differential $d$ on $\lasu_X$  is completely determined by the following:
\begin{itemize}
\item[(1)]
 The  $0$-simplices of $X$ are Maurer-Cartan elements.
\item[(2)]
 The linear part of $d$ is the desuspension of the differential in $N_*(X)$.
\item[(3)]
If $j\colon Y\subset X$ is a subsimplicial set, then
$
\lasu_j=\libc\bigl(s^{-1}N_*(j)\bigr)$.
\end{itemize}

The category $\catcdgl$ inherits by transfer from $\catss$ a {\em cofibrantly generated model structure}  \cite[Chapter 8]{bufemutan0}, for which the functors in (\ref{pair}) become a Quillen pair with the classical model structure on $\catss$. Quasi-isomorphisms  are weak equivalences and fibrations are surjections only for cdgl morphisms between connected cdgl's. In particular,
a quasi-isomorphism of connected cdgl's of the form
$$
(\libc(V),d)\stackrel{\simeq}{\longrightarrow} L
$$
makes of  $(\libc(V),d)$ a cofibrant replacement of $L$ and we say that it is a {\em  (Lie) model of $L$}. If $d$ has no linear term we say that $(\libc(V),d)$ is {\em minimal} and is unique up to cdgl isomorphism.

A fundamental object in this theory is the following \cite[Chapter 8.4]{bufemutan0}: let $X$ be a reduced simplicial set and let  $a$ be its only $0$-simplex.  The {\em minimal model of $X$} is, by definition, the  minimal model $(\libc(V),d)$ of  $\lasu_X^a$. In the same way, a {\em Lie model of $X$} is a model (not necessarily minimal) of $\lasu_X^a$. Observe that $\lasu_X^a$ is itself a Lie model although it is not minimal in general.

If $(\libc(V),d)$ is the minimal model of $X$ then $sV\cong \widetilde H_*(X;\bq)$. More generally, if $\varphi \colon (\libc(V),d)\to(\libc(W),d)$ is a Lie model of the map $f\colon X\to Y$ in $\catss_0$, then $H_*(f;\bq)$ is identified with the morphism $H(Q\varphi)\colon H(V,d_0)\to H(W,d_0)$ induced by $\varphi$ in the homology of the indecomposables of the models.

We finally remark, see \cite{fefuenmu0}, that the Quillen pair (\ref{pair2}) extends, up to homotopy, the classical functors of Quillen \cite{qui}. Moreover, for any cdgl $L$, its realization $\langle L\rangle$ is a strong deformation retract of $\mc_\bullet(L)$ which, in turn, extends the {\em Deligne groupoid} of $L$.

\section{A Lie characterization of the Bousfield-Kan $\bq$-completion}

This section is mostly devoted to the proof of Theorem \ref{main}. We begin with an observation of general nature, most of it well known, which will be often used. Recall that a Lie algebra is said to be {\em degreewise nilpotent} if for each $p\in\bz$ there exists an integer $n\ge1$ such that $L^n_p=0$.

\begin{proposition}\label{nilpo} For a connected complete Lie algebra $L$ the following are equivalent:

\begin{itemize}
\item[{\em(i)}] $L_0$ is a nilpotent group with the BCH product and the exponential action on $L$ is nilpotent.

    \item[{\em(ii)}] $L_0$ is a nilpotent Lie algebra and acts nilpotently on $L$.

    \item[{\em(iii)}] $L$ is degreewise nilpotent.
    \end{itemize}

    \end{proposition}

\begin{proof} The equivalence of {(ii)} and {(iii)} is an easy exercise.
On the other hand,  the well known Malcev equivalence between nilpotent Lie algebras and $\bq$-local nilpotent groups tells us precisely that $L_0$ is a nilpotent Lie algebra if and only if it is a nilpotent group with the BCH product. Also,
 if $L_0$ acts nilpotently on $L$ in the usual way, it does so via the exponential action as the successive commutators of this action raises the bracket length.

 Conversely, write $\varphi_x=e^{\ad_x}-\id$ for any $x\in L_0$. Consequently, the exponential action is nilpotent  if  there is an integer $n\ge 1$ such that
 $$
 \varphi_{x_1}\circ\dots\circ\varphi_{x_n}=0\quad\text{
for any $x_1,\dots,x_n\in L_0$.}
$$
 Note that
 $$
 \ad_x=\log(e^{\ad_x})=\log(\varphi_x+\id)=\sum_{k\ge 1}\frac{(-1)^{k+1}}{k}\varphi_x^k,
 $$
 so that, in
 $$
 \ad_{x_1}\circ\dots\circ\ad_{x_n}=\Bigl(\sum_{k\ge 1}\frac{(-1)^{k+1}}{k}\varphi_{x_1}^k\Bigr)\circ\dots\circ \Bigl(\sum_{k\ge 1}\frac{(-1)^{k+1}}{k}\varphi_{x_n}^k\Bigr),
 $$
 any summand vanishes. Hence, $
 \ad_{x_1}\circ\dots\circ\ad_{x_n}=0$ for any  $x_1,\dots,x_n\in L_0$ which implies that $L_0$ acts nilpotently on $L$ in the usual way.
\end{proof}

\begin{definition}\label{homolog} A connected cdgl $L$ is {\em homologically nilpotent} if $H(L)$ satisfies any of the equivalent properties of the previous result.
\end{definition}

 For the remainder of the section we fix $L=(\libc(V),d)$, a connected minimal cdgl. 

\subsection{A model of $\langle L/L^n\rangle$ is $L/L^n$}

In this subsection  we prove:

\begin{proposition}\label{modelonil}
For each $n\ge 2$ the counit
$$
\varepsilon\colon \lasu^a_{\langle L/L^n\rangle}\stackrel{\simeq}{\longrightarrow} L/L^n
$$
is a quasi-isomorphism. In other words, the minimal model of ${\langle L/L^n\rangle}$ is that of $L/L^n$.
\end{proposition}

For it we need:

\begin{lemma}\label{lema}
 The model and realization functors induce equivalences
\begin{equation}\label{homo}
 \xymatrix{ \Ho\catss^{\mathcal Nf}_{0,\bq}& \Ho\catcdgl^{\mathcal Nf}_0 \ar@<1ex>[l]^(.47){\langle\,\cdot\,\rangle}
\ar@<1ex>[l];[]^(.50){\lasu^a}\\}
\end{equation}
between the homotopy categories of reduced, nilpotent, rational simplicial sets of finite type and  connected, homologically nilpotent cdgl's whose homology is of finite type.
\end{lemma}

\begin{proof}
First, note that for any  connected $cdgl$ $M$, its realization $\langle M\rangle$ is a $\bq$-local simplicial set in the Bousfield sense \cite[\S2]{bous}. That is,  any map $f\colon X\to Y$ inducing an isomorphism in rational homology induces a weak homotopy equivalence
$$
f^*\colon \map(Y,\langle M\rangle)\stackrel{\simeq}{\longrightarrow} \map(X,\langle M\rangle).
$$
 Indeed, choose $\varphi\colon A\stackrel{\simeq}{\longrightarrow} B$ a connected cdga model of $f$ which is necessarily a quasi-isomorphism. Then\footnote{The {\em complete tensor product} of a cdga $A$ and a cdgl $M$ is defined as $A\otimesc M=\varprojlim_{n} A\otimes M/F^n$ where $\{F_n\}_{n\ge 1}$ is the filtration for which $M$ is complete.}, $\varphi\,\otimesc\id_M\colon A\otimesc M\stackrel{\simeq}{\longrightarrow}B\otimesc M$ is also a quasi-isomorphism of connected cdgl's so that its realization
$$
\langle \varphi\,\otimesc\id_M\rangle \colon\langle  A\otimesc M\rangle \stackrel{\simeq}{\longrightarrow}\langle B\otimesc M\rangle
$$
is a weak homotopy equivalence which is weakly equivalent to $f^*$
  by \cite[Theorem~6.6]{ber}, see also \cite[Theorem~12.18]{bufemutan0}.

  Now, in view of (\ref{homoto}) and (\ref{action}) together with Proposition \ref{nilpo}, the realization  $\langle M\rangle$ of any  $M\in\catcdgl_0^{\mathcal Nf}$   is  a reduced, nilpotent simplicial set of finite type which, being $\bq$-local, is then rational.

  On the other hand, any reduced, nilpotent, rational simplicial set $X$ of finite type is weakly equivalent to $\langle \lasu^a_X\rangle$ by \cite[Theorem 11.14]{bufemutan0}. As a result, again by (\ref{homoto}) and (\ref{action}) and Proposition \ref{nilpo}, $\lasu^a_X\in \catcdgl_0^{\mathcal Nf}$ and  the unit $\eta$ of the adjunction (\ref{homo}) is  equivalent to the identity.

  Finally,  for any $M\in \catcdgl_0^{\mathcal Nf}$, and denoting by $\varepsilon$ the counit of (\ref{homo}), it follows by the appropriate triangle identity that $\langle\varepsilon(M)\rangle\circ\eta\langle M\rangle =\id_{\langle M\rangle}$. Thus, $\langle\varepsilon(M)\rangle$ is a weak equivalence which implies that $\varepsilon(M)\colon \lasu^a_{\langle M\rangle}\stackrel{\simeq}{\longrightarrow} M$ is a quasi-isomorphism and the lemma follows.
  \end{proof}

\begin{proof}[Proof of Proposition \ref{modelonil}]
 Let $\{L_{[i]}\}_{i\in I}$ be the family of finitely generated sub cdgl's of $L$  and write
$$
L=\cup_{i\in I} \,L_{[i]}=\varinjlim_{i\in I}\,L_{[i]}.
$$
Then, for each $n\ge 2$, $L_{[i]}/L_{[i]}^n$ is a finite type sub cdgl of $L/L^n$ and
$$
L/L^n= \cup_{i\in I}\,  L_{[i]}/L_{[i]}^n= \varinjlim_{i\in I}\, L_{[i]}/L_{[i]}^n.
$$
Moreover, $\langle L_{[i]}/L_{[i]}^n\rangle$ is a sub simplicial set of $\langle L/L^n\rangle$ and
$$
\langle L/L^n\rangle=\cup_{i\in I}\,\langle L_{[i]}/L_{[i]}^n\rangle =\varinjlim_{i\in I}\,\,\langle L_{[i]}/L_{[i]}^n\rangle .
$$
Indeed a $q$-simplex of $\langle L/L^n\rangle$ is a cdgl morphism $\lasu_q\to L/L^n$ whose image is clearly finitely generated so it factors through some $L_{[i]}/L_{[i]}^n$.
Therefore, since the model functor preserves inductive limits,
$$
\lasu_{\langle L/L^n \rangle}^a\cong \varinjlim_{i\in I}\,\lasu_{\langle L_{[i]}/L_{[i]}^n\rangle}^a .
$$
Finally, by Lemma \ref{lema}, the counit of (\ref{homo}) provides  quasi-isomorphisms
$$
\varepsilon_i\colon \lasu_{\langle L_{[i]}/L_{[i]}^n\rangle}^a\stackrel{\simeq}{\longrightarrow}L_{[i]}/L_{[i]}^n
$$
and thus
$$
\varepsilon=\varinjlim_{i\in I}\,\varepsilon_i\colon \lasu_{\langle L/L^n \rangle}^a\stackrel{\simeq}{\longrightarrow}L/L^n
$$
is also a quasi-isomorphism.
\end{proof}

\subsection{The minimal model of $L/L^n$ with a quadratic differential}\label{seccion}

In this subsection we assume that the differential in $L=(\libc(V),d)$ is quadratic, i.e., $dV\subset [V,V]$ and denote  by $\libc^m(V)\cong L^m/L^{m+1}$ the subspace of $L$ generated by elements  of bracket length $m$ so that $L=\prod_{m\ge1} \libc^m(V)$. Note that the upper degree on $L$ can also be obtained by declaring $V$ to be concentrated in upper degree $1$, i.e.,  $V=V^1$ and extend it bracketwise to $\libc(V)$. Since $d \,\libc^m(V)\subset \libc^{m+1}(V)$ it follows that the homology of $L$ inherits an {\em upper}  grading for which
$$H(L)\cong\prod_{m\ge 1}H^m(L).
$$

We fix $n\ge 2$ and  determine the minimal model of $L/L^n$  following the approach in \cite[\S2]{5authors} for its counterpart in the ``commutative setting''.

The bigrading in $L$ trivially induces a bigrading in $L/L^n$ for which  the projection
$$
\varrho \colon L\longrightarrow L/L^n
$$
is also bigraded. Observe that, in the upper grading, $H^m(\varrho)$ is an isomorphism if $m\le n-2$ and is injective if $m=n-1$.

Next, choose a homogeneous bigraded basis $\{[a_i]\}$ of $\coker H^{n-1}(\varrho)$, with representatives $a_i\in (L/L^n)^{n-1}$ and define the bigraded vector space $U^{n-1}$ with basis $\{u_i\}$, concentrated in upper degree $n-1$, where each $u_i$ has  the same lower degree as $a_i$. Extend the bigrading on $U^{n-1}$  bracketwise to $\libc(V\oplus U^{n-1})$  declaring, as before,  $V=V^1$. Set $dU^{n-1}=0$, and  extend $\varrho$ to a bigraded morphism
$$
\varphi\colon (\libc(V\oplus U^{n-1}),d)\longrightarrow L/L^n,\qquad\varphi(u_i)=a_i,
$$
so that $H^{\le n-1}(\varphi)$ is an isomorphism.
Note also that  $(L/L^n)^n=0$ and therefore
$$
H^n(\varphi)\colon H^n(\libc(V\oplus U^{n-1}),d)\longrightarrow 0
$$
is the trivial map.
Choose then a homogeneous bigraded basis $\{[b_j]\}$ of $H^n(\libc(V\oplus U^{n-1}),d)$ and define a bigraded vector space $W_{[1]}^{n-1}$ with basis $\{w_j\}$, concentrated in upper degree $n-1$, where each $w_j$ has lower degree $|b_j|+1$. Extend this bigrading bracketwise to  $\libc(V\oplus U^{n-1}\oplus W^{n-1}_{[1]})$, set $dw_j=b_j$ and extend  $\varphi$ to
$$
\varphi\colon (\libc(V\oplus U^{n-1}\oplus W^{n-1}_{[1]}),d)\longrightarrow L/L^n
$$
so that it vanishes in $W^{n-1}_{[1]}$.

Recursively,  for each $r\ge 2$, assume that we have constructed  
$$
\varphi\colon (\libc(V\oplus U^{n-1}\oplus W^{n-1}_{<{[r]}}),d)\to L/L^n.
$$
We then define 
$$
W^{n-1}_{[r]}\cong  H^n(\libc(V\oplus U^{n-1}\oplus W^{n-1}_{<{[r]}}),d), 
$$
and extend both $d$ and $\varphi$ to $W^{n-1}_{[r]}$ as before.
Finally, define
\begin{equation}\label{descomposi}
W^{n-1}=\oplus_{r\ge 1}W^{n-1}_{[r]},\qquad Z^{n-1}=U^{n-1}\oplus W^{n-1}
\end{equation}
 and
$$
\varphi\colon (\libc(V\oplus Z^{n-1}),d)\longrightarrow L/L^n.
$$
By construction,  $d$ is of upper degree 1,
$$
H^n(\libc(V\oplus Z^{n-1}),d)=0,
$$
and $H^k(\varphi)$ is an isomorphism if $k\le n$ and vanishes for $k\ge n$.

Assume that we have built $Z^{n-1},\dots,Z^{m-1}$ and a morphism
$$
\varphi\colon (\libc(V\oplus Z^{n-1}\oplus\dots\oplus Z^{m-1}),d)\to L/L^n
$$
such that $d$ is of upper degree $1$,
$$
H^m(\libc(V\oplus Z^{n-1}\oplus\dots\oplus Z^{m-1}),d)=0,
$$
and $H^k(\varphi)$ is an isomorphism if $k\le m$ and vanishes for $k\ge m$.

 Repeating the above recursive procedure, for each $r\ge 1$, assume that we have constructed a morphism
$$
\varphi\colon (\libc(V\oplus Z^{n-1}\oplus\cdots\oplus Z^{m-1}\oplus W^m_{<{[r]}}),d)\to L/L^n.
$$
We then define
$$
W^m_{[r]}\cong H^{m+1}\bigl(\libc(V\oplus Z^{n-1}\oplus\cdots\oplus Z^{m-1}\oplus W^m_{<{[r]}}),d\bigr),
$$
concentrated in upper degree $m$, and extend both $d$ and $\varphi$ to $W^m_{[r]}$ as above. Finally, define
$$
Z^{m}=\oplus_{r\ge 1}W^{m}_{[r]}
$$
 and
$$
\varphi\colon (\libc(V\oplus Z^{n-1}\oplus\dots\oplus Z^{m}),d)\longrightarrow L/L^n.
$$
By construction, the differential $d$ is of upper degree 1, and $H^k(\varphi)$ is an isomorphism if $k\le m+1$ and vanishes for $k\ge m+1$.

We have thus constructed
$$
Z=\oplus_{m\ge n-1}Z^m
$$
and a quasi-isomorphism
\begin{equation}\label{minimalmodel}
\varphi\colon (\libc(V\oplus Z),d)\stackrel{\simeq}{\longrightarrow}L/L^n
\end{equation}
in which $d$ is of bidegree $(1,-1)$ with respect to the upper and lower gradings.

Furthermore, this model of $L/L^n$ satisfies the following properties:

\begin{proposition}\label{remark1}
\begin{itemize}
\item[{\em(i)}]$Z$ is concentrated in positive (homological degrees), i.e., $Z=Z_{\ge 1}$.

\item[{\em(ii)}] For any $p\ge 1$ and any $m\ge n-1$,
$$
dZ^m_p\subset \libc(V\oplus Z^{\le m}_{<p}).
$$
\item[{\em(iii)}]
$H^{\ge n }(\libc(V\oplus Z),d)=0$.

\item[{\em(iv)}] The differential in $\libc(V\oplus Z)$ is decomposable so that this is the minimal model of $L/L^n$. Moreover,
$$
d\,Z^{n-1}\subset [V,Z^{n-1}]\oplus \libc^{n}(V).
$$
\end{itemize}
\end{proposition}

\begin{proof}
{(i)} Since $\libc(V)$ is connected the map
$$
H_0(\varrho)\colon H_0(L)\longrightarrow H_0(L/L^n)
$$
is surjective for any $n$. In particular $\coker H_0^{n-1}(\varrho)=0$ and thus $U^{n-1}_0=0$.
On the other hand, by construction, the differential on the remaining  summands of $Z^{n-1}$ and on $Z^m$ for  $m\ge n$   is nontrivial. Hence, these spaces are concentrated in positive homological degrees.

\smallskip

Statements (ii) and (iii) follow directly from the  construction of $(\libc(V\oplus Z),d)$.

\smallskip

{(iv)} Recall that $dU^{n-1}=0$. On the other hand, since $d$ is of upper degree $1$, and $V$ is concentrated in upper degree $1$ there are no indecomposable cycles of upper degree $n$ in $(\libc(V\oplus U^{n-1}),d)$. Hence $dw$ is decomposable for any $w\in W^{n-1}_{[1]}$. For the same reason there are no indecomposable cycles of upper degree $n$ in $(\libc(V\oplus U^{n-1}\oplus W^{n-1}_{<[r]}),d)$ and again $dw$ is decomposable for any $w\in W^{n-1}_{[r]}$. This shows that $dz$ is decomposable for any $z\in Z^{n-1}$.

As for $m>n-1$,  taking into account (ii) and that $d$ is of upper degree $1$, there are no indecomposable elements of upper degree $m+1$ in $\libc(V\oplus Z^{\le m-1})$. It follows that $d$ must also be decomposable on $Z^m$.

In particular, as $d$ increases the upper degree by $1$, the subspace $dZ^{n-1}$ consists  of decomposable elements of $\libc(V\oplus Z^{n-1})^n$. Moreover, since $[Z^{n-1},Z^{n-1}]\subset Z^{2n-2}$, and $2n-2>n$ for $n\ge 3$, upper degree considerations show that no such term can contribute to upper degree $n$. It follows that
$$
d\,Z^{n-1}\subset [V,Z^{n-1}]\oplus \libc^{n}(V).
$$
For $n=2$ the same conclusion holds by construction, since $dZ^1$ has no component in $[Z^1,Z^1]$.
\end{proof}

Another specific feature of this model, to which we pay special attention and which will be extremely useful in what follows, is that elements of $Z$ with low homological degree are also concentrated in low upper degree. We emphasize that the specific upper bound in the following statement is not important. What really matters is the existence of such a bound depending only on
 $n$ and on the homological degree.

\begin{proposition}\label{lakey}
For any $p\ge 1$,
$$
Z_p\subset Z^{\le k_{n,p}}\quad\text{with}\quad k_{n,p}= (p+1)(n-2)+1.
$$
\end{proposition}
\begin{proof}
Consider the adjoint pair of  classical Quillen functors \cite[Appendix~B]{qui}
$$
\xymatrix{
{\rm \bf cdgc} \ar@<0.75ex>[r]^-{\mathscr L} &\mathbf{dgl} \ar@<0.75ex>[l]^(0.49){{\mathscr C} }
}
$$
between the categories of (counital and coaugmented) cocommutative differential graded coalgebras and differential graded Lie algebras. Recall that
$${\mathscr L}(C)= (\mathbb L(s^{-1}\overline{C}),d),$$
where $\overline{C}$ denotes the coaugmentation ideal, and $d=d_1+d_2$ whith
$$
d_1(s^{-1}c) = -s^{-1}dc\quad\text{and}\quad
d_2(s^{-1}c) = \frac{1}{2} \sum_i (-1)^{\vert a_i\vert} [s^{-1}a_i, s^{-1}b_i] $$
for $\overline{\Delta}c = \sum_i a_i\otimes b_i$. On the other hand,
$${\mathscr C}(L) = (\land (sL),d),$$
where again $d=d_1+d_2$ with
$$d_1(sv) = -sdv\quad\text{and}\quad d_2(sv\land sw)= (-1)^{\vert v \vert+1} s[v,w].$$

\medskip
Now, to avoid  cumbersome notation, write $M=L/L^n$ and recall from \cite[Prop.~2.3]{bufemutan0} that the counit of the above adjunction
$$
\gamma\colon\quil\quic(M)=(\lib(s^{-1}\Lambda^+ s M),d)\stackrel{\simeq}{\longrightarrow} M,
$$
which is the unique dgl morphism extending the projection
$$
s^{-1}\Lambda^+ s M\longrightarrow s^{-1}\Lambda^+ s M/(s^{-1}\Lambda^{\ge 2} s M)\cong M,
$$
is a quasi-isomorphism.
Write $u(x)$ for the upper degree of an element $x\in M$ and extend it to any generator
$$
\alpha=s^{-1}(sx_1\wedge\dots\wedge sx_r)\in s^{-1}\Lambda^+ s M
$$
by
$$
u(\alpha)= \sum_{i=1}^r u(x_i)-r+1.
$$
Observe that, since $M=M_{\ge 0}$,
$$
|\alpha|=\sum_{i=1}^r |x_i|+r-1\ge r-1.
$$
In addition,   since every element of $M$ necessarily has upper degree less than $n$
(because $M^n=0$), we have
$$
u(\alpha)=\sum_{i=1}^r u(x_i)-r+1\le (n-2)r+1.
$$
Hence, if $|\alpha|=p$, we conclude that $r\le p+1$ and therefore,
\begin{equation}\label{grado}
u(\alpha)\le (n-2)(p+1)+1.
\end{equation}
Now extend the upper degree on $s^{-1}\Lambda^+ s M$  bracketwise to
$\lib(s^{-1}\Lambda^+ s M)
$,
where a straightforward check shows  that $d$ is of degree 1 with respect to the upper grading. Consequently,  $H(\lib(s^{-1}\Lambda^+ s M),d)$ is bigraded. Furthermore,
 $\gamma$ is a morphism of bigraded dgl's.

 We next observe that the completion of the quasi-isomorphism $\gamma$,
$$
\widehat\gamma\colon (\libc(s^{-1}\Lambda^+ s M),d)\stackrel{\simeq}{\longrightarrow} \widehat M=M
$$
is also a quasi-isomorphism.  On the one hand, $H(\widehat\gamma)$ is clearly surjective. On the other hand, any element in $\libc(s^{-1}\Lambda^+ s M)$ of bidegree $(p,q)$, with respect to the homological and upper gradings, is necessarily an element of the uncomplete dgl $\lib(s^{-1}\Lambda^+ s M)$. Since $H(\libc(s^{-1}\Lambda^+ s M),d)$ is bigraded, this readily implies that $H(\widehat\gamma)$ is also injective. We conclude that $\widehat\gamma$ is a Lie model of $L/L^n$.

Next, following \cite[Prop.~3.18]{bufemutan0}, decompose $s^{-1}\Lambda^+ s M$ as
$$
s^{-1}\Lambda^+ s M\cong B\oplus dB\oplus A
$$
so that
$$
(\libc(A),d) \subset (\libc(s^{-1}\Lambda^+ s M),d)
$$
is minimal. Thus, the restriction of $\widehat\gamma$ to
$$
\widehat\gamma\colon (\libc(A),d)\stackrel{\simeq}{\longrightarrow} M
$$
is the minimal Lie model of $L/L^n$. By uniqueness of minimal models and (iv) of Proposition \ref{remark1}, it follows that
$$
A\cong  V\oplus Z
$$
and therefore, up to bigraded isomorphisms, $\widehat\gamma$ is precisely the morphism $\varphi$ in (\ref{minimalmodel}). Furthermore, note that both maps restrict to the identity on  $V$.

Finally, a direct inspection shows that the upper grading on $Z$, viewed as a subspace of $s^{-1}\Lambda^+ s M$, coincides with the original. Therefore, equation (\ref{grado}) completes the proof.
\end{proof}

\subsection{A  model of $L/L^n$}

We now attack the general case and consider $L=(\libc(V),d)$ a minimal cdgl in which the differential $d$ is not necessarily quadratic so that it decomposes  as
\begin{equation}\label{deri}
d=\sum_{m\ge 1}d_m\quad\text{with}\quad d_mV\subset \libc^{m+1}(V).
\end{equation}
In particular $(\libc(V),d_1)$ is a minimal cdgl with quadratic differential. Hence, for a fixed $n\ge 2$, we have, as in (\ref{minimalmodel}), a bigraded quasi-isomorphism
\begin{equation}\label{quasi2}
\varphi\colon (\libc(V\oplus Z),d_1)\stackrel{\simeq}{\longrightarrow}(L/L^n,d_1).
\end{equation}
In what follows, we denote by $\libc(V\oplus Z)^m$ the subspace of $\libc(V\oplus Z)$ spanned by elements of upper degree exactly $m$. Recall that the upper grading is given by extending bracketwise the grading $V=V^1$ and  $Z=\oplus_{r\ge n-1}Z^r$. This should not be confused with the usual space $\libc^m(V\oplus Z)$ which is generated by Lie brackets of length $m$, regardless  of the upper degree of each of the elements involved.

\begin{proposition}\label{perturba} There is a differential $d$ in $\libc(V\oplus Z)$  extending that of $\libc(V)$ and $d_1$ in $\libc(V\oplus Z)$ such that
$$
d=\sum_{m\ge 1}d_m, \quad d_mZ^r\subset \libc(V\oplus Z)^{r+m},\quad m\ge 1,\quad r\ge n-1,
$$
and
$$
\varphi\colon (\libc(V\oplus Z),d)\stackrel{\simeq}{\longrightarrow}L/L^n
$$
is a quasi-isomorphism.
\end{proposition}

\begin{proof}
We follow a perturbation argument to build $d$ in $\libc(V\oplus Z)$ of the required form
$$
d=\sum_{m\ge 1}d_m
$$
satisfying:
\begin{itemize}
\item[$\scriptstyle\bullet$] $d$ extends the original differential in $L=\libc(V)$. In other words,  $d_m$ extend the corresponding derivation in (\ref{deri}) for any $m\ge 1$.
  \item[$\scriptstyle\bullet$] For each $m\ge 1$,  $d_m$ is a derivation.
  \item[$\scriptstyle\bullet$] For each $m\ge1$ and each $r\ge n-1$,
  $$d_mZ^r\subset \libc(V\oplus Z)^{r+m}.
  $$
      \item[$\scriptstyle\bullet$] $d^2=0$. That is,
          $$\sum_{i+j=m}d_id_j=0,\qquad  m\ge 1.
          $$
      \end{itemize}
We set $d_1$ to be the differential in (\ref{quasi2}). By construction, see \S\ref{seccion}, $d_1$ increases the upper length by 1 so that
$$
d_1Z^r\subset \libc(V\oplus Z)^{r+1}, \quad r\ge n-1.
$$

\smallskip

We now define $d_m$, for $m\ge 2$, on $Z^r$, for $r\ge n-1$ by induction on $m$, on $r$, and on the lower (homological) degree of each $Z^{r}$ which is always positive by  (i) of Proposition \ref{remark1}.

\smallskip

From now on, we fix a homogeneous bigraded basis of $Z$, and whenever we refer to an element of $Z_p^r$,
 we implicitly mean one of the elements of this basis.

\smallskip

We begin by defining $d_2$ on $Z^{n-1}$ and fix an element $z\in Z^{n-1}_1$. By (iv) of Proposition \ref{remark1},
 $$
 d_1z\in \libc(V)^n_0\oplus [Z^{n-1},V]_0=\libc(V)^n_0,
 $$
 since $[Z^{n-1},V]_0$ vanishes for degree reasons. Hence
  $d_2d_1z$ is a well defined element in $\libc(V)^{n+2}$ and it is a $d_1$-cycle since $d_1d_2d_1z=-d_2d_1^2z=0$.
  However, by (iii) of Proposition \ref{remark1}, any $d_1$-cycle in $\libc(V)^{n+2}$ is a $d_1$-boundary in $\libc(V\oplus Z)$. That is, $d_2d_1z=d_1\Phi$ for some $\Phi\in \libc(V\oplus Z)^{n+1}$. We define $d_2z=-\Phi$ so the identity
  $$
  (d_2d_1+d_1d_2)(z)=0
   $$
   holds\footnote{Caution: $\Phi$ might well be an indecomposable element of $Z^{n-1}$. Hence the differential of the statement is not minimal in general.}.

  \smallskip

Assume now that $d_2$ has been defined on $Z^{n-1}_{<p}$ for $p\ge 1$ and let $z\in Z^{n-1}_p$. Again, by (iv) of Proposition \ref{remark1},
$
d_1z$ lies in $ \libc(V)^n_{p-1}\oplus [Z^{n-1},V]_{p-1}
$
where $d_2$ is already defined. Using the same argument as before  $d_2d_1 z$ is a $d_1$-cycle in $\libc(V\oplus Z)^{n+2}$ and thus it is a $d_1$-boundary in $\libc(V\oplus Z)$: $d_2d_1 z=d_1\Phi$ with $\Phi\in \libc(V\oplus Z)^{n+1}$. Define $d_2z=-\Phi$. This completes the induction step and $d_2$ is now defined on $Z^{n-1}$ with the required properties.

  \smallskip

Now suppose that $d_2,\dots,d_{m-1}$ have been defined on $Z^{n-1}$ and we proceed to define $d_m$ in $Z^{n-1}$ again by induction on the lower (homological) degree.

Let $z\in Z^{n-1}_1$.
As noted above $d_1z\in \libc(V)^{n}$ so that $d_md_1z$ is already defined and so is $$
(d_2d_{m-1}+\cdots+d_md_1)(z).
$$
This element is again  a $d_1$-cycle as a short computation shows, taking into account that  $\sum_{i+j=m-1}d_id_j=0$ by  the inductive hypothesis.
Thus, once again  by (iii) of Proposition \ref{remark1},
 $$
 (d_2d_{m-1}+\cdots+d_md_1)(z)=d_1\Psi
  $$
 for some   $\Psi\in \libc(V\oplus Z)^{n+m-1}$.  Define $d_mz=-\Psi$ so that $$
 \sum_{i+j=m}d_id_j(z)=0.
 $$
 Assume now that $d_m$ has been defined in $Z^{n-1}_{<p}$ for some $p\ge 1$ and let $z\in Z^{n-1}_p$. Again, by (iv) of Proposition \ref{remark1},
$$
d_1z\in \libc(V)^n_{p-1}\oplus [Z^{n-1},V]_{p-1}
$$
where $d_m$ is already defined. Using the same argument as before  $$
(d_2d_{m-1}+\cdots+d_md_1)(z)=d_1\Psi
$$
for some  $\Psi\in \libc(V\oplus Z)^{n+m+1}$. Define $d_mz=-\Psi$. This completes the induction and $d$ is defined in $Z^{n-1}$ with  the required properties.

\smallskip

 For $r\ge n$, the same procedure applies to define $d$ in $Z^r$, where we now invoke (ii) of Proposition \ref{remark1},
 $$
 d_1Z^r_p\subset \libc(V\oplus Z^{\le r}_{<p})
 $$
 and $d_md_1Z^r_p$ is then inductively defined for any $m$ and any $p\ge 1$.

 \smallskip

We have completed the definition of $d$ as claimed, but it remains to check that, for any $z\in Z$, the series
\begin{equation}\label{welldefined}
\sum_{m\ge 1}d_mz
\end{equation}
defines an element in $\libc(V\oplus Z)$.
 To this end note first that, for any $k$ and $\ell$,
 $$
 \libc(V\oplus Z^{\le \ell})^k\subset \libc^{\ge k/\ell}(V\oplus Z).
 $$
Now,  by Proposition \ref{lakey} and for any $p\ge 1$, we have $Z_{ p}\subset Z^{\le k_{n,p}}$ with $k_{n,p}=(p+1)(n-2)+1$, for any $p\ge 1$.
It follows  that, for any $z\in Z_p^r$,
$$
d_mz\subset \libc(V\oplus Z_{<p})^{r+m}\subset \libc(V\oplus Z^{\le  k_{n,p}})^{r+m}\subset \libc^{\ge (r+m)/k_{n,p}} (V\oplus Z).
$$
Consequently, only finitely many summands $d_mz$ contribute to any fixed bracket length and hence (\ref{welldefined}) is well defined.

 Note also that, by construction, $d$ commutes with $\varphi$ so that
$$
\varphi\colon (\libc(V\oplus Z),d) \longrightarrow L/L^n
$$
is indeed a cdgl morphism.

Finally, consider in the domain and codomain of this map the filtrations given by $\libc(V\oplus Z)^{\ge q}$ and $(L/L^n)^{\ge q}$ for $q\ge 1$. Then, $\varphi$ respects the filtration and the induced map  between the $E_1$-pages of the corresponding spectral sequences is precisely the quasi-isomorphism
$$
\varphi\colon (\libc(V\oplus Z),d_1)\stackrel{\simeq}{\longrightarrow}(L/L^n,d_1).
$$
\end{proof}

\begin{rem}\label{nominimal} As observed in the proof of Proposition \ref{perturba}, the model $(\libc(V\oplus Z),d)$ need  not be minimal as $d$ is not decomposable in general: only its restriction to $\libc(V)$, which coincides with the original differential, and the component $d_1$ are necessarily decomposable. The fact that we  write $d=\sum_{m\ge 1} d_m$ should not confuse the reader: each $d_m$ increases the upper degree by $m$, but not necessarily the usual bracket length.
\end{rem}

Nevertheless,  despite the possible non-minimality of $(\libc(V\oplus Z),d)$, the decomposable part of the differential, denoted by $\partial$, remains sufficiently well controlled for our purposes and satisfies the following property.

\begin{lemma}\label{aux1} $
 H(V\oplus Z,\partial)\cong V\oplus R
$
for some subspace $R\subset Z$.
\end{lemma}

\begin{proof} As observed in the above remark the differential in $V$ is always decomposable and, by construction, no element of $V$ lies in the image of  $\partial$.  Therefore, for  the decomposition
$$
V\oplus Z\cong C\oplus B\oplus \partial B,
$$
in which $C\cong  H(V\oplus Z,\partial)$,
it follows that
$
C\cong V\oplus R
$ for some  $R\subset Z$.
\end{proof}

\subsection{The proof}

We complete here the proof of Theorem \ref{main}. To start, fix a reduced simplicial set $X$, let $a$ be the only $0$-simplex, and consider $X\to  \langle\lasu^a_X\rangle$ the unit of the adjunction (\ref{pair2}). For the rest of the section denote by  $L=(\libc(V),d)$  the minimal model of $X$. Then $X\to  \langle\lasu^a_X\rangle$ is homotopic to a map of the form
$$
X\longrightarrow \langle L\rangle\cong\langle \varprojlim_n L/L^n\rangle\cong \varprojlim_n\langle L/L^n\rangle,
$$
which is induced by a commutative diagram

\begin{equation}\label{tower}
\xymatrix{
&&&X\ar[dll]_{f_1}\ar[dl]^(.40){f_2}\ar[d]^{f_3}\ar[drr]^{f_n}& &\\
& {*}&\langle L/L^2\rangle\ar[l]^(.54){p_1} &\langle L/L^3\rangle\ar[l]^(.47){p_2}&\cdots\ar[l]&\langle L/L^n\rangle\ar[l]&\cdots\ar[l]^(.38){p_n}
 }
\end{equation}
The core of the proof is to show that this is a {\em $\bq$-tower} for $X$ in the sense of \cite[Chapter III, \S6.1]{bouskan}. That is, provided  $p_n$ is a fibration for each $n\ge 1$, the diagram induces a pro-isomorphism for each $i\ge 1$,
$$
\{H_i(X;\bq)\}_{n\ge1}\cong\{H_i(\langle L/L^n\rangle;\bq)\}_{n\ge1}.
$$
In other words, see \cite[Chapter III, \S2.1]{bouskan} or the original source \cite{arma}, for each $i\ge1$ and each $n\ge 1$, there exists an integer $m\ge n$ and a map $H_i(\langle L/L^m\rangle;\bq)\to H_i(X;\bq)$ such that the following commutes:
\begin{equation}\label{diaim}
\xymatrix{
H_i(X;\bq)\ar@{=}[rr]\ar[d]_{H_i(f_n)}&&H_i(X;\bq)\ar[d]^{H_i(f_m)}\\
H_i(\langle L/L^n\rangle;\bq)&&H_i(\langle L/L^m\rangle;\bq).\ar[llu]\ar[ll]^{H_i(p_{n}{\scriptscriptstyle\circ}\dots{\scriptscriptstyle\circ} p_{m-1})}}
\end{equation}

\smallskip

Note that, in view of (\ref{homoto}), (\ref{action}) and Proposition \ref{nilpo},  the simplicial sets $\langle L/L^n \rangle$ are  nilpotent for $n\ge 1$. Moreover, see the proof of Lemma \ref{lema}, these are all rational spaces so that $\langle L/L^n \rangle\simeq \langle L/L^n \rangle^\acento_\bq$. Hence, by Lemmas 6.3 or 6.4 of  \cite[Chapter III, \S6]{bouskan}, if (\ref{tower}) were a $\bq$-tower, then $\langle L\rangle\simeq \varprojlim_n\langle L/L^n\rangle$ would be weakly homotopy  equivalent to $X^\acento_\bq$ and Theorem \ref{main} would follow.

\medskip

Begin by noting that each projection $q_n\colon L/L^{n+1}\to L/L^n$ is a cdgl fibration since it is surjective. Therefore, $p_n=\langle q_n\rangle$ is a Kan fibration of simplicial sets.

Next, by Proposition \ref{modelonil}, and for each $n\ge 2$, a Lie model of $f_n$ is given by the projection $L\to L/L^n$. Furthermore, this projection can be factored as
$$
\xymatrix{
L\,\ar@{^(->}[r]\ar[rd]&(\libc\bigl(V\oplus Z(n)\bigr),d)\ar[d]^{\varphi_n}_\simeq\\
&L/L^n}
$$
where $\varphi_n$ and $d=\sum_{m\ge 1}d_m$ are as in Proposition \ref{perturba}.
\begin{proposition}\label{propove} For each $n\ge 2$ there is a morphism

$$
 \rho^n\colon (\libc\bigl(V\oplus Z(n+1)\bigr),d)\longrightarrow(\libc\bigl(V\oplus Z(n)\bigr),d)
 $$
which fixes $\libc(V)$, it makes the  diagram
$$
\xymatrix{
(\libc\bigl(V\oplus Z(n+1)\bigr),d)\ar[d]_{\varphi_{n+1}}^\simeq\ar[r]^(.53){\rho^n}&(\libc\bigl(V\oplus Z(n)\bigr),d)\ar[d]^{\varphi_{n}}_\simeq\\
L/L^{n+1}\ar[r]_{q_n}&L/L^{n},
}
$$
commutative and can be written as
$$
\rho^n=\sum_{m\ge 0}\rho^n_m,\quad\text{with}\quad \rho^n_m\bigl(Z^r(n+1)\bigr)\subset \libc\bigl(V\oplus Z(n)\bigr)^{r+m}\quad\text{ for all $r\ge n$.}
$$
\end{proposition}
\begin{proof}
For simplicity in the notation fix $n\ge2$ and replace $Z(n+1)$, $Z(n)$, $q_n$, $\varphi_n$, $\varphi_{n+1}$ and  $\rho^n$ by $Y$, $Z$, $q$, $\varphi$, $\psi$ and $\rho$ respectively, so that the morphism to be constructed has the form
$$
\rho\colon (\libc(V\oplus Y),d)\longrightarrow (\libc(V\oplus Z),d)
$$
and it is required to fit in the diagram
$$
\xymatrix{
(\libc\bigl(V\oplus Y\bigr),d)\ar[d]_{\psi}^\simeq\ar[r]^{\rho}&(\libc\bigl(V\oplus Z\bigr),d)\ar[d]^{\varphi}_\simeq\\
L/L^{n+1}\ar[r]^{q}&L/L^{n}.
}
$$
Note that the map $\rho$, with the prescribed decomposition
\begin{equation}\label{ro}
\rho=\sum_{m\ge 0}\rho_m,\quad\text{with}\quad \rho_m(Y^r)\subset \libc\bigl(V\oplus Z\bigr)^{r+m}\quad\text{ for all $r\ge n$.}
\end{equation}
 commutes with differentials if and only if
\begin{equation}\label{perturbf}
\sum_{i=0}^m\rho_i\,\,d_{m+1-i}=\sum_{i=0}^md_{m+1-i}\,\,\rho_i,\quad m\ge 0.
\end{equation}
 Recall from (\ref{descomposi}) that $Y^n=U^n\oplus W^n$ where the differential vanishes on $U^n$. Then, set $\rho(U^n)=0$. To define $\rho$ on $W^n\oplus Y^{>n}$ we follow a perturbation argument completely analogous to the one used in the proof of Proposition \ref{perturba} and therefore only sketch its construction.

Let $\rho_0,\dots,\rho_{k-1}$ be constructed on $Y$ so that (\ref{ro}) and (\ref{perturbf}) hold for  $r\ge n$ and $0\le m\le k-1$. Assume also $\rho_k$ defined on $Y^{\le s}_{<p}$ for some $s\ge n$ and some $p\ge 1$ so that  (\ref{ro}) and (\ref{perturbf}) are satisfied for $n\le r\le s$ and $m=k$. Then, for any $y\in Y^s_p$, a short computation shows that
$$
\bigl(\sum_{i=0}^k\rho_i\,\,d_{m+1-i}-\sum_{i=0}^{k-1}d_{m+1-i}\,\,\rho_i)(y)
$$
is a $d_1$-cycle in $\libc(V\oplus Z)^{s+k+1}$ which is sent to zero by the quasi-isomorphism $\varphi$. Thus, it is  $d_1\Phi$ with $\Phi\in \libc(V\oplus Z)^{s+k}$. Define $\rho(y)=\Phi$  and the induction is complete.

To check that
$$
\rho=\sum_{m\ge 0}\rho_m
$$
 is well defined, we follow the same argument as in the proof of Proposition \ref{perturba}. Recall from Proposition \ref{lakey} that for any $p\ge 1$,
$$
Z_{\le p}\subset Z^{\ge k_{n,p}}\quad \text{with}\quad k_{n,p}=(p+1)(n-2)+1.
$$
  Therefore, for any $p\ge 1$, $m\ge 0$ and $r\ge n$,
$$\rho_m(Y_p^r)\subset \libc(V\oplus Z_{\le p})^{r+m} \subset  \mathbb L^{\geq (r+m)/k_{n,p}} (V\oplus Z).$$
It follows  that, for any $y\in Y_p^r$ only finitely many summands $\rho_mz$ contribute to a fixed bracket length, and hence $\rho(y)=\sum_{m\ge 0}\rho_m(y)$ is well defined.

\end{proof}

\begin{lemma}\label{keyrem}
Given $p\ge 1$ and $n\ge 2$, there is always an integer $m>n$ such that the composition
$$
\rho=\rho^{n}{\scriptstyle\circ}\dots {\scriptstyle\circ}\,\rho^{m-1}\colon  (\libc\bigl(V\oplus Z(m)\bigr),d)\longrightarrow(\libc\bigl(V\oplus Z(n)\bigr),d)
$$
is decomposable on $Z(m)_p$.
\end{lemma}

\begin{proof}
By Proposition \ref{lakey}, $Z(n)_{\le p}\subset Z(n)^{\le k_{n,p}}$. Choose $m> k_{n,p}$. Then, in view of Proposition \ref{propove}, for any $z\in Z(m)_p$, and taking into account that $Z(m)$ is concentrated in upper degrees greater than or equal to $m-1$, we obtain
$$
\rho(z)\in \libc\bigl(V\oplus Z(n)\bigr)^{\ge m-1}_p\subset \libc\bigl(V\oplus Z(n)^{\le k_{n,p}}\bigr)^{>k_{n,p}}
$$
and hence $\rho(z)$ must be decomposable.
\end{proof}
Next, in view of Proposition \ref{modelonil}, and for each $m>n\ge 2$, each triangle
$$
\xymatrix{
&X\ar[dl]_{f_n}\ar[dr]^{f_{m}}\\
\langle L/L^n\rangle&&\langle L/L^{m}\rangle\ar[ll]^{p_{n}{\scriptscriptstyle\circ}\dots{\scriptscriptstyle\circ}p_{m-1}}
}
$$
is modeled by
$$
\xymatrix{
&\quad\quad\,\,(\libc(V),d)\quad\quad\,\,\ar@{_(->}[dl]\ar@{^(->}[dr]&\\
(\libc\bigl(V\oplus Z(n)\bigr),d)&&(\libc\bigl(V\oplus Z(m)\bigr),d)\ar[ll]^{\rho}
}
$$
where, as in Lemma \ref{keyrem}, we write
$$
\rho=\rho^{n}{\scriptstyle\circ}\dots {\scriptstyle\circ}\,\rho^{m-1}.
$$
Consequently, the diagram
$$
\xymatrix{
&s^{-1}\widetilde H_*(X;\bq)\ar[dl]_{s^{-1}\widetilde H_*(f_n)\quad}\ar[dr]^{\quad s^{-1}\widetilde H_*(f_{m})}\\
s^{-1}\widetilde H_*(\langle L/L^n\rangle;\bq)&&s^{-1}\widetilde H_*(\langle L/L^{m}\rangle;\bq)\ar[ll]^{s^{-1}\widetilde H_*(p_{n}{\scriptscriptstyle\circ}\dots{\scriptscriptstyle\circ}p_{m-1})}
}
$$
is identified with the corresponding one on homology of indecomposables
$$
\xymatrix{
&\,\, V\,\ar[dl]\ar[dr]\\
H(V\oplus Z(n),\partial)&&H(V\oplus Z(m),\partial)\ar[ll]^{H\bigl(Q(\rho)\bigr)}
}
$$
where $Q(\rho)\colon (V\oplus Z(m),\partial)\to (V\oplus Z(n),\partial)$ denotes the map induced by $\rho$ on indecomposables. By Lemma \ref{aux1}, this diagram is of the form
$$
\xymatrix{
&\,\, V\,\ar@{_(->}[dl]\ar@{^(->}[dr]\\
V\oplus R(n)&&V\oplus R(m)\ar[ll]^{H\bigl(Q(\rho)\bigr)}
}
$$
for certain subspaces  $R(n)\subset Z(n)$ and $R(m)\subset Z(m)$. Furthermore, $H\bigl(Q(\rho)\bigr)$ necessarily fixes $V$.

Now, for a fixed homological degree $p\ge1$, choose $m>n$ as in Lemma \ref{keyrem} so that $\rho\bigl(Z(m)_p\bigr)$, and thus $\rho \bigl(R(m)_p\bigr)$, is decomposable. In particular $H_p\bigl(Q(\rho)\bigr)\bigl(R(m)\bigr)=0$ and therefore,
 the projection $\bigl(V\oplus R(m)\bigr)_p\to V_p$ fits in the commutative diagram
 $$
\xymatrix{
&\,\,V_p\ar@{_(->}[dl]\\
\bigl(V\oplus R(n)\bigr)_p&&\bigl(V\oplus R(m)\bigr)_p.\ar[ll]^{H\bigl(Q(\rho)\bigr)}\ar[ul]
}
$$
In other words, for each $i\ge 1$ and $n\ge 1$, we have found $m>n$ together with a commutative diagram as in (\ref{diaim})
$$
\xymatrix{
H_i(X;\bq)\ar@{=}[rr]\ar[d]_{H_i(f_n)}&&H_i(X;\bq)\ar[d]^{H_i(f_{m})}\\
H_i(\langle L/L^n\rangle;\bq)&&H_i(\langle L/L^{m}\rangle;\bq).\ar[llu]\ar[ll]^{\, H_i(p_{n}{\scriptscriptstyle\circ}\dots{\scriptscriptstyle\circ} p_{m-1})\,}}
$$
Therefore, (\ref{tower}) is a $\bq$-tower of $X$ and this finishes the proof of Theorem \ref{main}.
\subsection{On nilpotent spaces}

When restricting to nilpotent spaces Theorem \ref{main} implies the following:

\begin{theorem}\label{conse1}
 The model and realization functors induce equivalences
 $$
 \xymatrix{ \Ho\catss^{\mathcal N}_{0,\bq}& \Ho\catcdgl^{\mathcal N}_0 \ar@<1ex>[l]^(.47){\langle\,\cdot\,\rangle}
\ar@<1ex>[l];[]^(.50){\lasu^a}\\}
$$
between the homotopy category of reduced nilpotent rational simplicial sets  and  homologically nilpotent connected cdgl's.
\end{theorem}

\begin{proof} Recall that $\bq$-completion coincides with rationalization in the homotopy category of connected nilpotent spaces and the latter is idempotent. Hence, by Theorem \ref{main}, the unit $X\stackrel{\simeq}{\longrightarrow}\langle \lasu^a_X\rangle$ is a weak homotopy equivalence for any rational nilpotent reduced simplicial set. In view of (\ref{homoto}), (\ref{action}) and Proposition \ref{nilpo}, this also shows that $\lasu^a_X\in \catcdgl^{\mathcal N}_0$.

On the other hand, recall from the proof of Lemma \ref{lema}, that for any  connected $cdgl$ $L$, its realization $\langle L\rangle$ is a $\bq$-local simplicial set. In particular, if $L\in \catcdgl^{\mathcal N}_0$, then $\langle L\rangle$ is nilpotent and thus rational.
To finish, as in the proof of Lemma \ref{lema}, the suitable triangle identity provides that the counit $\lasu^a_{\langle L\rangle}\stackrel{\simeq}{\longrightarrow} L $ is a quasi-isomorphism for any $L\in \catcdgl^{\mathcal N}_0$.
\end{proof}

\subsection{A final remark}
Since the Bousfield $\bq$-homology localization \cite{bous} is closely related to the $\bq$-completion as it coincides with it for $\bq$-good spaces, we have deemed it appropriate to include the following observation. In $\catss$ consider the model category structure introduced by Bousfield in \cite[\S10]{bous} in which weak equivalences (we$_\bq$) are rational homology isomorphisms and cofibrations are, like in the classical structure, injections. Hence,
if we denote by we, fib and cof the weak equivalences, fibrations and cofibrations in the usual structure, it follows that we $\subset$ we$_\bq$, cof $ = $ cof$_\bq$ and fib$_\bq$ $\subset$ fib. Moreover,  trivial fibrations are the same in both structures.

Note then that the realization and model functors also constitute a Quillen pair when $\catss$ is considered with the new structure. Indeed, $\lasu$ trivially preserves cofibrations. On the other hand, $\varphi$ is a fibration in $\catcdgl$ if and only if $\langle \varphi\rangle$ is a Kan fibration, whose domain and codomain are $\bq$-local.  By \cite[Proposition~3.3.16(1)]{hirsch} it is also a fibration in the Bousfield model structure.

\bigskip

\noindent{\bf Acknowledgment.} The authors are grateful to the referee for a careful reading of the manuscript and for numerous suggestions and corrections that significantly improved both the exposition and the mathematics.

\smallskip

The authors also thank Prof. Joost Nuiten for a helpful discussion.

\bigskip
\noindent {\sc Institut de Math\'ematiques et Physique, Universit\'e Catholique de Louvain, Chemin du Cyclotron 2,
1348 Louvain-la-Neuve,
         Belgique}.

\noindent\texttt{yves.felix@uclouvain.be}

\medskip

\noindent{\sc  Departamento de  Geometr\'{\i}a y Topolog\'{\i}a, Universidad de Sevilla, C/ Tarifa s/n, 41012 Sevilla, Spain.}

\noindent
\texttt{mario.fuentes.rumi@gmail.com}

\medskip

\noindent{\sc Departamento de \'Algebra, Geometr\'{\i}a y Topolog\'{\i}a, Facultad de Ciencias, Universidad de M\'alaga, Blvr. Louis Pasteur 31, 29010 M\'alaga, Spain.}

\noindent
\texttt{aniceto@uma.es}

\end{document}